\documentclass [12pt,leqno]{article}
\usepackage[leqno]{amsmath}
\usepackage{paralist,amsthm,amsfonts,amssymb,graphics, topcapt}%
\usepackage{t1enc, pst-poly, pstricks,color,appendix}
\usepackage[hypertex]{hyperref}
\usepackage{hyperref}
\divide
\oddsidemargin by 2 
\makeatletter
\newcommand{\address}[2][]{%
  \ifx\@add@ress\@undefined\gdef\@add@ress{\par\par\bigskip}\AtEndDocument{\@add@ress}\fi
  \g@addto@macro\@add@ress{\bigskip\noindent{\small\scshape%
      \ifx#1\empty\else{\bfseries Address of #1:}\ \fi#2}\par\par}}

\newcommand{\email}[1]{\noindent{\normalfont{\itshape E-mail:\/} \texttt{#1}}\par\par}
\renewenvironment{abstract}{\small\quotation\noindent
  {\bfseries \abstractname}}{\endquotation \par}
\newcommand{\footnotetextplain}[1]{\begingroup\def\@thefnmark{}%
  \long\def\@makefntext##1{\parindent 0pt\noindent ##1}\@footnotetext{#1}
  \endgroup}

\newcommand{\AMSsubjclass}[2]{\footnotetextplain{2000
   \emph{Mathematics Subject Classification:} Primary #1, Secondary #2.}}

\newcommand{\keywords}[1]{\footnotetextplain{\emph{Key words and phrases:} #1.}}

\makeatletter \xdef\qedbuit{\qed}

\newcommand{\TeoremaAmbFinalMarcat}[1]{%
  \expandafter\gdef\csname end#1\endcsname{\qedbuit\@endtheorem}}

\newtheorem{theo}{Theorem}[section]

\theoremstyle{definition}
 \TeoremaAmbFinalMarcat{defi}
 \TeoremaAmbFinalMarcat{ex}
\newtheorem{rem}[theo]{Remark} \TeoremaAmbFinalMarcat{rem}
 \TeoremaAmbFinalMarcat{hrem}
 \TeoremaAmbFinalMarcat{conv}
 \TeoremaAmbFinalMarcat{prob}
\newtheorem{conj}[theo]{Conjecture} \TeoremaAmbFinalMarcat{conj}
 \TeoremaAmbFinalMarcat{conj}

\newenvironment{proclama}[1]{
                \par\vspace{\topsep}\noindent{\bf #1}
                \begin{em}}
                {\end{em}\par\vspace{\topsep}}

\newcommand{\start}[2]{\begin{#1}\label{#2}}

\newcommand{\theoc}[1]{Theorem~\ref{#1}}
\newcommand{\propc}[1]{Proposition~\ref{#1}}

\newcommand{\lemc}[1]{Lemma~\ref{#1}}

\newcommand{\figc}[1]{Figure~\ref{#1}}

\newcommand{\defc}[1]{Definition~\ref{#1}}

\def\@enum@{\list{\csname label\@enumctr\endcsname}%
           {\usecounter{\@enumctr}\def\makelabel##1{\hss\llap{##1}}
           \itemsep=2pt\parsep=0pt\topsep=3pt plus 1pt minus 1 pt}}

\newenvironment{numlist}{\enumerate[(1)]}{\endenumerate}

\def\map#1#2#3{\mbox{${#1}\colon {#2} \longrightarrow {#3}$}}
\def\Smap#1#2{\mbox{${#1}\colon{#2} \longrightarrow {#2}$}}

\def\u{\mathsf{U}}
\def\v{\mathsf{V}}
\def\w{\mathsf{W}}
\def\x{\mathsf{X}}
\def\y{\mathsf{Y}}

\def\p{\mathsf{P}}
\def\q{\mathsf{Q}}

\def\UU{\mathcal{U}}
\def\VV{\mathcal{V}}
\def\WW{\mathcal{W}}

\def\V{\mathbb{V}}
\newcommand{\vvv}{\pi^\ast}
\def\W{\mathbb{W}}
\def\bb{\mathbf{LP}}

\def\ov#1{\overline{#1}}

\def\ovv{\overline{v}}

\def\wh#1{\widehat{#1}}
\def\wv{\widehat{\v}}

\def\id{\mathop\mathrm{Id}}
\def\sign{\mathrm{sign}}

\def\cob{\delta}
\def\cib{\Delta}

\def\om{\omega}
\def\sss{\mathsf{s}}

\newcommand{\aq}{\ensuremath{\mathbb{A}}_q}

\newcommand{\Z}{\ensuremath{\mathbb{Z}}}

\newcommand{\m}{\mathsf{M}}

\def\lu{\ell_\mathsf{U}}
\def\lv{\ell_\mathsf{V}}
\def\lw{\ell_\mathsf{W}}

\newcommand{\sequence}{%
{\rput{18}(1.4,1.4){\PstPolygon[unit=1.3,PolyNbSides=10] }
\rput{20}(3.9,1.4){\PstPolygon[unit=1.3,PolyNbSides=10] }
\rput{20}(6.3,1.4){\PstPolygon[unit=1.3,PolyNbSides=10]}
\rput{20}(8.8,1.4){\PstPolygon[unit=1.3,PolyNbSides=10]}
\rput{20}(11.3,1.4){\PstPolygon[unit=1.3,PolyNbSides=10]} } }

\newcommand{\xx}{
 \rput(-0.5,1.3){$\w$}
 \rput(-0.5,0.3){$\u$}
 \rput(0,1.5){
 \rput(1.25,0){$\v_k$}
 \rput(3,0){$\v_k$}
 \rput(5,0){$\v_k$}
\rput(6.75,0){$\v_{\scriptscriptstyle k,h}$}}
\psline(0,2)(7.5,2)
\psline[linecolor=gray,linewidth=3pt](0,2)(0,1)
\psline(2,1)(2,2)
\psline(4,1)(4,2)
\psline[linewidth=3pt](6,1)(6,2)
\psline[linecolor=gray,linewidth=3pt](7.5,1)(7.5,2)
\psline(0,1)(7.5,1)

 \rput(1.3,0){
  \psline[linecolor=gray,linewidth=3pt](0.1,1)(0.1,0)
 \rput(1,0.3){$\v_i$}
 \rput(3,0.3){$\v_i$}
 \rput(5,0.3){$\v_i$}
  \rput(7.2,0.3){$\v_{\scriptscriptstyle i\!,\!j}$}
  \psline[linecolor=gray,linewidth=3pt](7.5,0)(7.5,1)%
\psline(7.5,1)(0.2,1)
\psline(2,0)(2,1)
\psline(4,0)(4,1)
\psline[linewidth=3pt](6,0)(6,1)
\psline(0.1,0)(7.5,0)
}

\psdot(0.2,0.7)
\psline{->}(0.2,1.2)(0.2,0.8)

\psdot[dotstyle=o](0.2,1.2)
\rput(0.25,1.45){$v_{k}$}

\psdot(5.6,0.7)
\psline{->}(5.6,0.8)(5.6,1.2)
\psdot(5.6,1.2)

\rput(5.6,1.4){$v_h$}
\psline{->}(7.5,0.8)(5.7,0.8)

\psdot[dotstyle=o](7.6,1.2)
\rput(7.8,0.4){$u_{\scriptscriptstyle \lu\!-\!r}$}

\psline{->}(7.6,1.1)(7.6,0.8)
\rput(6.6,0.6){${\scriptstyle \lv}$}

\psdot(7.6,0.7)
\rput(7.7,1.45){$v_k$}

  \psline[linecolor=gray]{[-}(0.1,0.4)(0.4,0.4) \rput(0.6,0.4){$r$}\psline[linecolor=gray]{-]}(0.7,0.4)(1.2,0.4)
}

\title{Algebraic characterization of simple closed curves via Turaev's cobracket}

\author{Moira Chas and Fabiana Krongold}
\address{\footnotesize Department of Mathematics, Stony Brook University, NY 11794, USA  \hfill\break \email{moira@math.sunysb.edu}}

\address{\footnotesize
Departamento de Ciencias Exactas, Ciclo B\'asico Com\'un,\hfill\break
Universidad de Buenos Aires
Buenos Aires, Argentina \hfill\break \email{fkrong@dm.uba.ar} }

\date{\today}

\begin{document}
\maketitle
\begin{abstract}  The vector space $\V$ generated by the conjugacy classes in the fundamental group  of an orientable surface has a natural Lie cobracket $\map{\delta}{\V}{\V\times \V}$. For negatively curved surfaces, $\delta$ can be computed from a geodesic representative as a sum over  transversal self-intersection points. In particular $\delta$ is zero for any power of an embedded
simple closed curve. Denote by Turaev(k)  the statement that $\delta(x^k) = 0$ if and only if the nonpower conjugacy class $x$ is represented by an embedded curve.
Computer implementation of the cobracket delta unearthed counterexamples to Turaev(1) on every surface with negative Euler characteristic except the pair of pants.
Computer search have verified Turaev(2) for hundreds of millions of the shortest classes. In this paper we prove Turaev(k) for $k=3,4,5,\dots$ for surfaces with boundary.
Turaev himself introduced the cobracket in the 80's and wondered about the relation with embedded curves, in particular asking if Turaev (1) might be true.

We give an application of our result to the curve complex. We show that a permutation of the set of free homotopy classes that commutes with the cobracket and the power operation is induced by an element of the mapping class group.

\end{abstract}
\AMSsubjclass{57M99}{17B62}
%
%
%


\section{Introduction}

In \cite{T}, Turaev defined a Lie coalgebra structure in the vector space generated by free homotopy classes of non-trivial closed curves on an orientable surface. The cobracket of a class is a sum over self-intersection points of representative of this class. Therefore, the cobracket of a simple closed curve is zero. It is not hard to see that the cobracket of a power of a simple closed curve is zero. Thus, he asked whether the converse is true.

In \cite{chas} it was shown that, if one excludes trivial cases and the pair of pants, the answer of this question is negative. The main goal of our work here is to show that, for surfaces with non-empty boundary, one may reformulate Turaev's question so it has a positive answer. In other words, Turaev's cobracket yields a characterization of simple closed curves. More precisely, we prove:

\begin{proclama}{Main Theorem} Let $\VV$ be a free homotopy class of curves on an oriented surface with boundary.
 Then  $\VV$ contains a power of a simple curve if and only if the Turaev cobracket of $\VV^3$ is zero. Moreover, if $\VV$ is nonpower and $p$ is an integer larger than three then the number of terms of the cobracket of $\VV^p$ (counted with multiplicity) equals $2p^2$ times the self-intersection number of $\VV$.\end{proclama}

Our main tool is the presentation of Turaev's cobracket in \cite{chas}. Using this presentation, we list the terms of the cobracket of powers of a cyclic word and we can show that if the power is large enough, these terms do not cancel.

We also give an application of our result to the curve complex.

This paper is structured as follows: We start by describing in Subsection \ref{Notation} the vector space of cyclic reduced words  and  combinatorial cobracket. In Subsection \ref{not conjugate}, using techniques similar to those of  in \cite{CK}, we prove that certain pairs of cyclic
words cannot be conjugate (\propc{main}).
This result distinguishes the conjugacy classes of various sets of linear words. In Section~\ref{turaev}, we use Proposition~\ref{main} to prove our Main Theorem. In Section~\ref{ccc} we give an application of our Main Theorem to the curve complex.
For completeness, in the appendix, we include the following know results: In Appendix~\ref{curves}, we recall the definition of Turaev's cobracket.
 In Appendix \ref{review c}, some results of \cite{chas} we use in this work are stated. Finally, in Appendix \ref{cc} we review the definition of the curve complex.

\section{The Turaev Lie coalgebra on the vector space generated by cyclic reduced words}

\subsection{Basic definitions and Notation}\label{Notation}
For each positive integer $q$,  a \emph{$q$-alphabet $\aq$}, or, briefly,
an \emph{alphabet}, is the set of  $2q$ symbols, called \emph{letters} $\{a_1,a_2,
\dots,a_q,\ov{a}_1, \ov{a}_2, \dots, \overline{a}_q\}$, endowed with
a fixed linear order.  For each letter $v$,
$\ov{\ov{v}}=v$.

A \emph{linear word in $\aq$}  is a  finite sequence of
symbols $v_0v_1\dots v_{n-1}$ such that $v_i$ belongs to $\aq$ for
each $i \in \{0,1,\dots,n-1\}.$  The empty word is a linear word with zero letters. Let $\v=v_0v_1\dots v_{n-1}$ be a linear word. By
definition, $\ov{\v}=\ov{v}_{n-1} \ov{v}_{n-2}\dots \ov{v}_0$.   The
linear word $\v$ is \emph{freely reduced} if $v_i \ne \ov {v}_{i+1}$
for each $i \in \{0,1,\dots,n-1\}$. If $\v$ is freely reduced and
$v_{n-1}\ne \ov{v}_0$ then $\v$ is \emph{cyclically reduced}.
Consider the equivalence relation on the set of linear words,
generated by the pairs of the form $(\v,\w)$ such that $\v$ is a
cyclic permutation of $\w$ or $\v=\w v \overline{v}$ where $v$ is a
letter in $\aq$. The equivalence classes  under this equivalence
relation  are called \emph{cyclic words}. (Observe that these are
the conjugacy classes of the free group generated by $a_1, a_2,
\dots, a_{q}$). If $\v$
is a (not necessarily reduced) linear word, we denote the
equivalence class of $\v$ by $\widehat{\v}$. A linear representative of $\widehat{\v}$ is a cyclically
reduced linear word in $\wv$.
A cyclic word $\widehat{\v}$ is  \emph{nonpower} if it does not contain a proper power of another word.
A linear subword of ${\widehat\v}$ is a linear subword of a linear representative of ${\widehat\v}$.
The \emph{length of $\v$}, denoted by $\ell_\v$ or $\ell(\v)$  is  the number of symbols $\v$ contains. The \emph{length of  the cyclic word $\widehat{\v }$} is  the length of a  linear representative of $\widehat{\v }$. Thus, if $\v$ is cyclically
reduced then the length of $\widehat{\v }$ equals $\ell_\v$.

 When dealing with letters denoting linear
words $\v=v_0v_1\dots v_{n-1}$ , subindices of
letters will be considered mod the length of $\v$, that is $n$.



From now on,   fix $\v=v_0v_1\dots v_{n-1}$, a cyclically reduced linear word of positive length $n$. The
set of linked pairs of the cyclic word $\wv$ is denoted by $\bb_1(\wv)$ (see appendix
\defc{linked}).

The next statement follows straightforwardly from  \defc{linked}.

\start{lem}{pair}
For each $({  \p, \q})$in  $\bb_1(\wv)$ there exist two integers  $i$ and $j$  in $\{0,1,\ldots,n-1\}$ such that the following holds:

\begin{numlist}
\item If $({  \p, \q})$ is a linked pair as in
\defc{linked} (1) or (2) then $\p=v_{i-r-1}v_{i-r}\dots v_{i}$ and
$\q=v_{j-r-1}v_{j-r}\dots v_{j}$, for some  integer $r \in
\{0,1,\ldots,n-2\}$.
\item If $({  \p, \q})$ is a linked pair as in \defc{linked} (3) then $\p=v_{i-r-1}v_{i-r}\dots v_{i}$ and
$\q=v_{j-1}v_j\dots v_{j+r}$,  for some  integer $r \in
\{1,2,\ldots, n-2\}$.
\end{numlist}
Moreover,  the map  $\phi: \bb_1(\wv)\to\{0,1,\ldots, n-1\} \times \{0,1,\ldots, n-1\}$ defined by $({
\p, \q})\mapsto (i,j)$ is injective.
\end{lem}


Since $\v$ is fixed, we can (and will)  identify each linked pair $(\p,\q)$ with the pair of integers $\phi(\p,\q)$ in \lemc{pair}.
Say that $(i,j)$ is a \emph{linked pair of type (1) or (2)} (respectively, \emph{(3)}) if $(i,j)=\phi(\p,\q)$ for some $(\p,\q)$ in $\bb_1(\wv)$ satisfying \defc{linked}(1) or (2) (respectively (3)).

\start{lem}{extremo}If $(i,j)$ is a linked pair of  $\wv$ then $v_{i}\ne v_{j}$ and
$\ov{v}_i\ne v_{j-1}$. In particular $i \ne j$.
\end{lem}

\begin{proof}
If $(i,j)$ is a linked pair of type (1) then  $(\overline{v}_{i-1},\overline{v}_{j-1},v_i,v_j)$  is positively
or negatively oriented. Hence $v_{i-1}\ne v_{j-1}$, $v_i\ne v_j$ and $\overline{v}_{j-1}\ne v_i$. If $(i,j)$ is a linked pair of type (2) or (3) the result follows straightforwardly using that $\v$ is cyclically reduced.
\end{proof}

For each pair $i$ and $j$  in $\{0,1,\dots, n-1\}$ denote
\[\v_{i,j}=\left\{\begin{array}{ll}v_{i} v_{i+1}\ldots v_{j-1} & \hbox{ if $i<j$, and}\\
                             v_{i}  v_{i+1}\ldots v_{n-1}v_{0}v_1\ldots v_{j-1} & \hbox{ if $i\ge j$.}\end{array}
            \right.\]
If $i=j$, $\v_i$ is used instead of $\v_{i,i}$.

The next result is a direct consequence of \defc{linked}.
\start{lem}{symmetry} \begin{numlist} \item If $(\p,\q)$ is a linked pair then $(\q,\p)$ is a linked pair.
However, $(i,j)$ and $(j,i)$ are both linked pairs  if an only if $(i,j)$ is a linked pair of type (1) or (2).
\item Let $(i,j)$ be a linked pair.  The linear word if $\v_{j,i}$ is cyclically reduced if an only if $(i,j)$ is a linked pair of type (1) or (2).
\item If $(i,j)$ is  a linked pair of type (3), then there exists  $r$ in $\{1,2,\dots,\ell_\v-1\}$ such that $r$ is the largest integer verifying  $v_{i-r}\dots v_{i-1}=\overline{v_{j}\dots v_{j+r-1}}$ and $(j+r,i-r)$ is a linked pair.
\end{numlist}
\end{lem}

The next result which establishes the exact relation between the linked pairs of a nonpower cyclic word $\widehat{\v}$ and its power $\widehat{\v}^p$, follows directly from \defc{linked}.

\start{lem}{linkedvp}
Let $\v$ be a cyclically reduced word of length $n$ and $p$ be a positive integer. Then,
$$\bb_1(\wv^p)=\Big\{\big(i+tn,j+sn\big): (i,j)\in \bb_1(\wv),  t,s\in\{0,1,2,\dots,p-1\}\Big\}.$$
\end{lem}

Denote by  $\V$  the vector space generated by the set of non-empty cyclic words in $\aq$.
The class of the empty word will be identified with the zero element in $\V$.

\start{prop}{delta} If $\v$ is cyclically reduced linear word and $p$ is a positive integer, then the map
$\map{\delta}{\V}{\V\otimes \V}$ of \defc{cobr} satisfies the following equality.
\begin{equation}{\label{eq delta}}
\delta(\widehat{\v}^p)= p \cdot\hspace{-5mm}\sum_{ \substack{(i,j)\in
                \bb_1(\wv) \\ 0\le s\le p-1}}\hspace{-5mm} \sign(i,j) \hspace{0.1cm} \wh{\Big(\v_i^s\v_{i,j}\Big)}\otimes
                 \wh{\Big(\v_j^{p-s-1}\v_{j,i} \Big)}.
\end{equation}
\end{prop}
\begin{proof}  Take $i, j \in \{0,1,\dots,n-1\}$ and $t, s \in \{0,1, \dots, p-1\}$. By \lemc{extremo}, $i \ne j$. It is not hard to see that
\[
\v^p_{i+tn,j+sn}=\left\{\begin{array}{ll}
				\v_i^{t-s}\ \v_{i,j} & \hbox{ if $i<j$ and $s \le t$,}\\
		         	\v_i^{p-t+s}\ \v_{i,j} & \hbox{ if $i<j$ and $s>t$,}\\	
                                 \v_i^{t-s-1}\ \v_{i,j} & \hbox{ if $i>j$ and $s<t$,}\\
                                  \v_i^{p-t+s-1}\ \v_{i,j} & \hbox{ if $i>j$ and $s \ge t$,}\\
                             	\end{array}
            \right.\]
The desired result follows from  \lemc{linkedvp}.
\end{proof}


\subsection{Certain words are not conjugate}\label{not conjugate}

In this subsection we prove combinatorially that certain pairs of cyclic words we construct out of  $\v$  cannot be
conjugate. Since the conjugacy classes of these words appear in the
terms of the cobracket $\delta(\widehat{\v}^p)$,
the fact that these words are not conjugate implies that these terms cannot cancel one
another.

 \start{prop}{main} Let $\v$  be a  cyclically reduced linear word  such that $\widehat{\v}$ is   nonpower and  $(i,j)$ and $(k,h)$ two distinct linked pairs of $\v$. If
$m$ is a positive integer greater than or equal to two, then, for
any positive integer $l$,
$$
 \widehat{\v_i^{m}\v}_{i,j}\not= \widehat{\v_k^{l}\v}_{k,h}
$$
\end{prop}

\begin{proof}

Let  $\u$  and $\w$ denote the linear words $\v_i^m\v_{i,j}$ and
$\v_k^l\v_{k,h}$ respectively.  By hypothesis, $\v$ is a cyclically
reduced linear word, therefore $\u$ and $\w$ are reduced. Moreover,
by \lemc{extremo}, $v_i\not=\overline{v}_{j-1}$ and
$v_k\not=\overline{v}_{h-1}$. Hence $\u$ and $\w$ are cyclically
reduced and
\begin{align} \label{d1}
m\cdot\lv< \lu <(m+1)\cdot \lv,\\
\nonumber l\cdot \lv< \lw <(l+1)\cdot \lv.
\end{align}

Set $\u=u_0u_1\ldots u_{\lu-1}$ and $\w=w_0w_1\ldots w_{\lw-1}$. The next equalities follow straightforwardly.
\begin{align}\label{p1}
 u_t&=u_{t+\lv} ,  \hbox{   for all $t$ in  $\{0, 1, \ldots , \lu -\lv-1\}$,}\\
\label{p2} w_t&=w_{t+\lv} ,  \hbox{   for all $t$ in $\{0, 1, \ldots , \lw -\lv-1\}$,}\\
\label{u} u_{\lu-\lv}&= v_{i+\lu-\lv}=v_j, \\
\label{w} w_{\lw-\lv}&=v_{k+\lw-\lv}=v_h.
\end{align}

We argue by contradiction: Assume that $\widehat\u= \widehat\w$.
Since the linear words $\u$ and $\w$ are linear representatives of
the same cyclic word, they have equal length, that is, $\lu=\lw$.
Furthermore, there exists $r$ in $\{0,1,\ldots,\lu-1\}$ such that
$\u=\w_{r}$. The proof is split into three cases:

\begin{enumerate}

\item $r=0$. Here, $\u=\w$. Hence the initial and final subwords of length
$\lv$ of $\u$ and $\w$ are equal i.e., $\v_i=\v_k$ and
 $\v_j=\v_h$. Distinct cyclic permutations of a nonpower cyclic reduced word are distinct, so $i=k$ and $j=h$, contradicting the hypothesis $(i,j)\ne (k,h)$.
\item    $0<r\le \lu-\lv$. By Equations (\ref{p2}) and (\ref{w}), since  $m$ is positive,
\[v_{k}=w_0=w_{\lw}=u_{\lw-r}= u_{\lw-r-\lv}=w_{\lw-\lv}=v_h,\]
which contradicts \lemc{extremo}. See \figc{Case 2}.

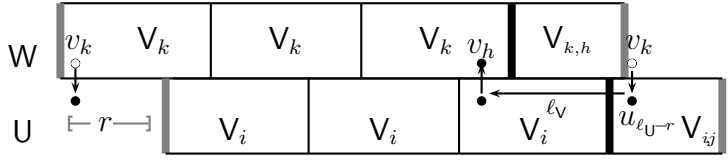
\begin{figure}[htbp]
 \begin{pspicture}(14,2)
 \rput(2,0){\xx}
   \end{pspicture}
   \caption{
   Case  $0<r\le \lu-\lv$ (here, $m=l=3$.)}
   \label{Case 2}
 \end{figure}

\item    $\lu-\lv< r< \lu$. Since $m \ge 2$, $\lv-1<\lu-\lv$. Hence $0\le \lu-r -1< \lu-\lv$.
By Equations (\ref{p1}) and (\ref{u}),
\[v_{j}=u_{\lu-\lv}=w_{\lw-\lv+r}=w_{\lw+r}=u_{\lu}=u_0=v_i,
\]
contradicting \lemc{extremo}.

\end{enumerate}
\end{proof}
\begin{rem}\label{ejemplo} The lower bound $2$ for $m$ in  \propc{main} is sharp.
 Indeed, consider the ordered alphabet $\{a,b,\overline{b},\overline{a}\}$ and the word $\v=abaabab$. The  linking  pairs $(0,1)$ and
$(5,6)$ are not equivalent. Nevertheless, the subword of $\v^2$ given by
$\v_0\v_{0,1}=abaababa$ is a cyclic permutation of $\v_5\v_{5,6}=ababaaba $. (Observe that the corresponding terms in the cobracket of $\v^2$,  $\widehat{\v_0\v}_{0,1}\otimes
\widehat{\v}_{1,0}$ and $\widehat{\v_5\v}_{5,6}\otimes \widehat{\v}_{6,5}$ are different because $\v_{1,0}=baabab$ and $\v_{6,5}= babaab$ which implies
$\widehat{\v}_{1,0} \ne \widehat{\v}_{6,5}=babaab$.)\end{rem}

 The \emph{Manhattan norm of an element $x$
of $\V \otimes \V$} denoted by $\m(x)$
 is the sum of the absolute values of the coefficients of the expression of
 $x$ in the basis of $\V\otimes \V$ consisting in the set of tensor products of pairs of cyclic words
 in $\aq$. Thus if $x=c_1 \VV_1\otimes\WW_1+c_2\VV_2\otimes \WW_2+\cdots+c_l \VV_l\otimes\WW_l$ where for each $i,j \in \{1,2,\dots,l\}$,
 $c_i \in \Z$, $\VV_i$ and $\WW_i$ are  cyclically reduced words and $\VV_i \ne \VV_j$ or $\WW_i \ne \WW_j$ when $i \ne j$ then $\m(x)=|c_1|+|c_2|+
 \cdots+|c_l|$.

\start{prop}{contar} Let $\v$  be a  cyclically reduced linear word and let $p$ an integer larger than three. Then the Manhattan norm of $\delta(\widehat{\v}^p)$ equals $p^2$ times
the number of elements in $\bb_1(\wv)$.\end{prop}

 \begin{proof} A simple argument shows that if the result holds for nonpower words, then it holds for all words. Therefore, we only need to prove it for nonpower words.
Assume that $\widehat{\v}$ is nonpower.
Let  $(i,j)$ be an element of $\bb_1(\wv)$. The terms of $\delta(\widehat{\v}^p)$ corresponding to $(i,j)$ in Equation (\ref{eq delta}) are
$\widehat{\v_i^m\v_{i,j}} \otimes \widehat{\v_j^{\hbox{\tiny $p\!\!-\!\!m\!\!-\!\!1$}}\v_{ji}}$ where $m \in \{0,1,2, \dots , p-1\}$.
Since the length of two of  the first factor of these elements differ by a non-zero multiple of $n$, these terms are all distinct.

We claim that if $(k,h) \in \bb_1(\wv)$,  $(k,h)\not=(i,j)$, and $l \in \{0,1,2,\dots, p-1\}$ then
$$\widehat{\v_i^m\v_{i,j}} \otimes \widehat{\v_j^{\hbox{\tiny $p\!\!-\!\!m\!\!-\!\!1$}}\v_{ji}} \ne \widehat{\v_k^l\v_{k,h}} \otimes \widehat{\v_h^{\hbox{\tiny $p\!\!-\!\!l\!\!-\!\!1$}}\v_{hk}}$$
Clearly, this claim implies the desired result.

To prove the claim observe that, since $p\ge 4$, either $m \ge 2$ or $p-m -1 \ge 2$. If $m\ge 2$, the result follows from by \propc{main}.  Hence, we can assume
 that $p-m -1 \ge 2$. If  $\widehat{\v_i^m\v_{i,j}}$ and $\widehat{\v_k^l\v_{k,h}}  $ are distinct the claim follows. Thus we can also assume that $\widehat{\v_i^m\v_{i,j}}=\widehat{\v_k^l\v_{k,h}}$  and therefore $m=l$. We complete the proof by showing that the assumption $\widehat{\v_j^{\hbox{\tiny $p\!\!-\!\!m\!\!-\!\!1$}}\v_{ji}} = \widehat{\v_h^{\hbox{\tiny $p\!\!-\!\!m\!\!-\!\!1$}}\v_{hk}}$
 leads to a contradiction.

Assume first that $(i,j)$ is  a linked pair of type (1) or (2). By \lemc{symmetry}(1)  $(j,i)$ is also a linked pair   a linked pair of type (1) or (2) and  $(h,k)\not=(j,i)$.
 By  \propc{main}, $\widehat{\v_j^{\hbox{\tiny $p\!\!-\!\!m\!\!-\!\!1$}}\v_{j,i}}  \ne \widehat{\v_h^{\hbox{\tiny $p\!\!-\!\!m\!\!-\!\!1$}}\v_{h,k}} $. Therefore, we can assume that  $(i,j)$ is a linked pair a linked pair of type (3). In this case, the length of $\widehat{\v_i^m\v_{i,j}} $ plus the length of  $\widehat{\v_j^{p-m-1}\v_{ji}}$ is strictly less than $p \cdot \ell_\v$. Furthermore, by  \lemc{symmetry}(3), there exists $r$ in $\{1,2,\dots,\ell_\v-1\}$ such that $v_{i-r}\dots v_{i-1}=\overline{v_{j}\dots v_{j+r-1}}$ and $(j+r,i-r)$ is a linked pair.  Hence,
 $
 \widehat{\v_j^{p-m-1}\v_{j,i}}=  \widehat{\v_{j+r}^{p-m-1}\v}_{j+r,i-r}
 $
 and the sum of the lengths of $\widehat{\v_i^m\v_{i,j}} $ and  $\widehat{\v_j^{p-m-1}\v_{ji}}$ equals $p \cdot \ell_\v-2r$. Therefore, the sum of the lengths of $\widehat{\v_k^l\v_{k,h}}$ and  $\widehat{\v_h^{\hbox{\tiny $p\!\!-\!\!l\!\!-\!\!1$}}\v_{hk}}$ is $p\cdot\ell_\v-2r$.  By \lemc{symmetry}(3),  $(k,h)$ is  a linked pair of type (3) and $(h+r,k-r)$ is linked. Hence,
$$\widehat{\v_j^{\hbox{\tiny $p\!\!-\!\!m\!\!-\!\!1$}}\v_{\hbox{\tiny
$\!\!j\!,\!i$}}}=\widehat{\v_{j+r}^{\hbox{\tiny $p\!\!-\!\!m\!\!-\!\!1$}}\v}_{\hbox{\tiny $j\!\!+\!\!r\!,\!i\!\!-\!\!r$}} =\widehat{\v_h^{\hbox{\tiny $p\!\!-\!\!l\!\!-\!\!1$}}\v_{h,k}}= \widehat{\v_{h+r}^{\hbox{\tiny $p\!\!-\!\!l\!\!-\!\!1$}}\v}_{\hbox{\tiny $h\!\!+\!\!r\!,\!k\!\!-\!\!r$}}.$$
By \propc{main}, $(j+r,i-r)=(h+r,k-r)$ and hence $(i,j)=(k,h)$ contradicting the hypothesis of the claim.

\end{proof}

\start{prop}{simple} For each cyclically reduced linear word $\v$,  if $\delta(\widehat{\v^3})=0$ then $\bb_1(\wv)$
is empty.
\end{prop}

\begin{proof}  Consider a cyclically reduced, linear, nonpower word $\w$ such that  $\v=\w^k$ for some positive integer $k$.
By \lemc{linkedvp}, $\bb_1(\wv)$ is empty if an only if \(\bb_1(\widehat{\w})\)  is empty. Therefore, if  $k>1$, the conclusion follows from \propc{contar}. Hence, we can assume that $\v$ is nonpower.
Assume  $\bb_1(\wv)$ is not empty.
 Let $(i,j)$ be an element in $\bb_1(\wv)$. By Equation (\ref{eq delta}), one of the terms of $\delta(\widehat{\v^3})$ is  $3\cdot \sign(i,j)\widehat{\v_i^2\v_{i,j}}\otimes \wh{\v_{j,i} }$. Denote this term by $T$. We will show that $T$  does not cancel with other terms of $\delta(\widehat{\v^3})$.

Observe that in  Equation (\ref{eq delta}) there are two other terms corresponding to $(i,j)$, namely,
$3\cdot \sign(i,j)\widehat{\v_{i,j}}\otimes \wh{\v_j^{2}\v_{j,i} }$ and $3\cdot \sign(i,j)\widehat{\v_i\v_{i,j}}\otimes \wh{\v_j\v_{j,i} }$.
Clearly, none of these terms cancel with $T$.

By \propc{main}, if $(k,h) \in \bb_1(\wv)$ and $(k,h) \ne (i,j)$ then
$\widehat{\v
_k^{l}\v_{\hbox{\tiny $\!\!k\!,\!h$}}}  \ne \widehat{\v_i^2\v_{i,j}}$  for any positive integer $l$. Implies that the terms corresponding to $(k,h)$ do not cancel with $T$. Thus the proof is complete.\end{proof}

\section{Reformulation of Turaev's conjecture}\label{turaev}

Let  $\Sigma$ be an orientable surface with non-empty boundary.  One knows that there is a bijection between the following three sets:
\begin{numlist}
\item the set $\vvv$ of non-empty,  cyclic reduced words on a set of free generators of the fundamental group of $\Sigma$ and their inverses,
\item the set of non-trivial, free homotopy classes of closed, oriented curves on $\Sigma$, and
\item the set of conjugacy classes of $\pi_1(\Sigma)$.
\end{numlist}
We will make use of this correspondence by identifying these three sets. Thus an element in $\vvv$ will sometimes be considered as a conjugacy class of  $\pi_1(\Sigma)$ and other times, a free homotopy class  of closed, oriented curves.

Given a cyclic reduced word $\VV \in \vvv$, define the
\emph{self-intersection number  $\sss(\VV)$} to be minimum number of transversal self-intersection points of a representatives of  $\VV$ which self-intersect only in  transverse double points.

By \theoc{iso}, the combinatorially defined map $\delta$ (see \defc{cobr}) and the Turaev cobracket $\Delta$ (see \defc{Turaevcobracket}) coincide. By \theoc{iso}, \propc{simple} and \propc{contar} we have the following result.

\begin{proclama}{Main Theorem} Let $\Sigma$ be an oriented surface with boundary and let $\VV$ in $\vvv$. Then  $\VV$ contains a power of a simple curve if and only if $\Delta(\VV^3)=0$. Moreover, if $\VV$ is nonpower and $p$ is an integer larger than three then the Manhattan norm of $\Delta (\VV^p)$ equals $2p^2$ times the self-intersection number of $\VV$. In symbols,
$
\m(\Delta (\VV^p))=2p^2\sss(\VV).
$
\end{proclama}

\section{The Curve complex and the Mapping Class Group}\label{ccc}

Say that an element $\VV \in \vvv$ is \emph{simple} if contains an embedded curve. The following lemma will be used in the proof of \propc{permutation}.

 \start{lem}{eight} If $\UU$ and $\VV$  be two disjoint simple elements of  $\vvv$
 then there exists $P$ in  $\Sigma$ and   $\u$ and $\v$ in $\pi_1(\Sigma,P)$ such that $\u \in \UU$, $\v \in \VV$ and the following holds.
 \begin{numlist}
        \item Either $\wh{\u.\v}$ or $\wh{\u.\overline{\v}}$ can be represented by a figure eight.
        \item If $\UU\not=\VV$ and $\wh{\u.\v}$ can be represented by a figure eight then $\wh{\u.\v}$ is nonpower and $\sss(\wh{\u.\v})=1$.
        \end{numlist}
\end{lem}

\begin{proof}  Let $\alpha$ and $\beta$ be disjoint closed simple curves that are representatives of $\UU$ and
$\VV$ respectively. Consider a connected component $C$ of $\Sigma\setminus(\alpha\cup\beta)$
including $\alpha\cup\beta$ in its boundary. Let $\gamma$ be a non self-intersecting  arc in $C$ from  a point $P$ in $\alpha$, to  a point $Q$ in $\beta$. Let $\u \in \pi_1(\Sigma,P)$ be the element  represented by $\alpha$ and let
$\v \in \pi_1(\Sigma,P)$ be the element represented by the curve that starts at $P$, runs along $\gamma$ to $Q$, then along $\beta$ and finally along $\gamma^{-1}$ from $Q$ to $P$.   It is not hard to see that $\u$ and $\v$ satisfy (1).

Now we prove (2). Since $\u.\v$ can be represented by a figure eight and $\UU$ and $\VV$ are non-trivial, $\sss(\wh{\u.\v})=1$. Consider an element $\WW \in \vvv$ such that $\wh{\u.\v}=\WW^p$,  for some positive integer $p$.  By definition of the cobracket,
$$
\Delta(\WW^p)=\Delta(\wh{\u.\v})=\pm(\UU \otimes \VV - \VV \otimes \UU) \ne 0.
$$
By \propc{power of simple}, $\WW$ is not simple, that is $\sss(\WW)>0$. By \cite[Theorem 3.9]{chas}, $ p^2\cdot s(\WW)\le
\sss(\WW^p)=\sss(\wh{\u.\v})=1$.
Then   $p=1$, and \(\wh{\u.\v}\) is
primitive.\end{proof}

\start{prop}{permutation} Let $\sigma$  be a  permutation on $\vvv$ which commutes
with the power operations (that is, for each $\VV$ in $\vvv$ and each $p$ in  $\Z$, $\sigma(\VV^p )=\sigma(\VV)^p$.) If the linear map induced by $\sigma$ on $\vvv$ commutes with the
Turaev's cobracket then
  \begin{numlist}
           \item the restriction of $\sigma$ to the subset of nonpower classes of a given self-intersection number is a permutation. In particular, $\sigma$ maps  simple classes into simple classes, and
        \item  the restriction of $\sigma$ to the subset of simple classes maps pairs of classes with disjoint representatives into pairs of classes with disjoint representatives.
    \end{numlist}
 Therefore, $\sigma$ induces an automorphism on the curve complex $C(\Sigma)$ (see Appendix~\ref{cc}).
\end{prop}

\begin{proof}
Since $\sigma$ is bijective and commutes with the power operations,  the restriction of $\sigma$ to the subset of nonpower classes is a permutation.

 If  $\UU \in \vvv$ is  nonpower, $\sigma(\UU)$ is nonpower. Since $\sigma$ is a bijection,  by the Main Theorem,
$$
2\cdot4^2\!\cdot\!\sss\big(\sigma(\UU)\big)=\m(\Delta(\sigma(\UU)^4)=
\m(\sigma(\Delta(\UU^4)))=
\m(\Delta(\UU^4))=2\cdot4^2\!\cdot\!\sss(\UU).
$$
 Therefore, $\sss(\UU)=\sss(\sigma(\UU))$. This completes the proof of (1).

Finally, consider two simple classes $\UU$ and $\VV$ in $\vvv$ that can be represented by disjoint closed curves.  If $\UU=\VV$, then the conclusion is immediate. Thus we can assume $\UU \ne \VV$. Let $\u$ and $\v$ in $\pi_1(\Sigma)$ be as   \lemc{eight} for $\UU$ and $\VV$. Suppose first that $\wh{\u\cdot\v}$ can be represented by a figure eight. By (1), $\sss(\sigma(\wh{\u\cdot\v}))=1$. Then there exists two elements $\x$ and $\y$ in $\pi_1(\Sigma)$ such that $\sigma(\wh{\u\cdot\v})=\wh{\x\cdot\y}$. Moreover, $\wh{\x}$ and $\wh{\y}$ have simple disjoint representatives. Hence,
$$\pm(\wh{\x}\otimes \wh{\y}-\wh{\y}\otimes \wh{\x})=\Delta(\sigma(\wh{\u\cdot \v}))=\sigma(\Delta(\wh{\u\cdot \v})) =\pm(\sigma( \wh{\u}) \otimes \sigma(\wh{\v})-\sigma(\wh{\v}) \otimes \sigma(\wh{\u}) )$$
Then  the set  $\{\wh{\x},\wh{\y}\}$ is equal to the set $\{\sigma(\wh{\u}),\sigma(\wh{\v})\}$. This implies that $\sigma(\wh{\u})$ and $\sigma(\wh{\v})$ are disjoint.

If $\wh{\u \cdot \overline{\v}}$ is represented by a figure eight, then by the arguments above we can show that $\sigma(\u)$ and $\sigma(\overline{\v})$ are disjoint. Since  $\sigma(\overline{\v})=\overline{\sigma(\v)}$, $\sigma(\u)$ and $\sigma(\v)$ are disjoint.
\end{proof}

We  denote by $\Sigma_{g,b}$ an oriented surface with genus $g$ and $b$ boundary components.

\start{theo}{comples} Let $\sigma$ be a permutation  on the set $\pi^\ast$ of free homotopy  classes of closed
curves on an oriented surface with boundary which is different from $\Sigma_{1,2}$, $\Sigma_{0,4}$ and $\Sigma_{1,1}$. Suppose the following:
\begin{numlist}
\item If $\sigma$ is extended linearly to the free $\Z$-module generated by $\vvv$, and to $\vvv \otimes \vvv$ (so $\sigma(\UU \otimes \VV)=\sigma(\UU)\otimes \sigma(\VV)$) then $\sigma$ preserves
the Turaev Lie cobracket. In symbols
$\Delta(\sigma(\UU))=\sigma(\Delta(\UU))$
for all $\UU \in \pi^\ast$.
\item For all $\UU \in \pi^\ast$, and all $p \in \Z$, $\sigma(\UU^p)=\sigma(\UU)^p$.
 \end{numlist}
Then $\sigma$ induces a permutation $\widetilde{\sigma}$ on the set of unoriented simple classes.
The permutation  $\widetilde{\sigma}$  is induced by a unique element of the mapping
class group. \end{theo}
\begin{proof}
By \propc{permutation}(1), since $\sigma(\overline{\UU})=\overline{\sigma(\UU)}$, $\sigma$ induces a permutation $\widetilde{\sigma}$ on the set of unoriented simple closed curves.

By \propc{permutation}(2), the restriction of $\widetilde{\sigma}$ to the set of unoriented simple closed curves preserves disjointness.

Thus, by \theoc{ikl}, $\widetilde{\sigma}$ is induced by an element of the Mapping Class group.
\end{proof}

This result ``supports'' Ivanov's statement in \cite{iv2}:

\begin{proclama}{Metaconjecture} ``Every object naturally associated with a surface $S$ and having a sufficiently rich
structure has $\mathrm{Mod}(S)$ as its group of automorphisms. Moreover, this can be proved by a reduction
theorem about the automorphisms of $\mathrm{C}(S)$.''
\end{proclama}

In this sense, the Turaev Lie cobracket combined with the power maps, have a ``sufficiently rich" structure.

\appendix

\section{Turaev's Lie coalgebra of curves on a surface}\label{curves}

Consider a vector space  $\W$. Define two linear maps $\Smap{\om}{\W\otimes \W\otimes \W}$ and $\Smap{\sss}{\W\otimes \W}$
by the formulae $\om(u \otimes v\otimes w) =w \otimes u \otimes v$ and $\sss(v\otimes
w)=w\otimes v$ for each triple of elements $u$, $v$ and $w$ in $\W$.
A {\em Lie cobracket} on $\W$ is a linear map $\map{\cib}{\W}{\W\otimes
\W}$ such that $\sss \circ\cib=-\cib$ ({\em co-skew symmetry}) and $(\id+\om+\om^2)(\id
\otimes \cib)\cib=0$ ({\em co-Jacobi identity}). If $\W$  is a vector space and $\cib$ is a cobracket on $\W$, $(\W,\cib)$ is a {\em Lie coalgebra}.

We recall the definition of Turaev's cobracket \cite{T}.

\start{defi}{Turaevcobracket} Let $\Sigma$ be an oriented surface. Choose a free set of generators for the fundamental group of $\Sigma$. The set of non-trivial free homotopy classes of curves on $\Sigma$ is in bijective correspondence with the set of cyclic  words (as defined in  Section \ref{Notation}) in the generators and their inverses. We will identify cyclic words with free homotopy classes.
Thus $\V$  is the vector space generated by the set of free homotopy classes of
non-trivial, oriented, closed curves on $\Sigma$.

Let $\gamma$ be a oriented, closed curve on $\Sigma$. Let $[\gamma]$ denote the free homotopy class of $\gamma$ if $\gamma$ is a non-contractible loop and $[\gamma]=0$ otherwise.
Consider a free homotopy class $[\beta]$ such that all self-intersection points of $\beta$ are  transversal double points.  Denote by $\mathcal{I}$ the set of self-intersection points of $\beta$. Each $P \in \mathcal{I}$ determines two loops based at $P$ so that $\beta$ can be obtained as the loop product of these two loops (forgetting the basepoint). Order these two loops so that the orientation given by the branch of the first, followed by the branch of the second equals the orientation of the surface. Denote the ordered pair of loops by   $(\beta_P^1,\beta_P^2)$.
The Turaev cobracket of $x$, $\Delta(x)$ is defined by the following formula:
\[
\Delta(x)=\sum_{P\in \mathcal{I}}
            \Big([\beta_P^1]\otimes [\beta_P^2]
            -[\beta_P^1]\otimes[\beta_P^2]\Big)
\]

\end{defi}
The next result follows straightforwardly from the definition of the cobracket.
\start{prop}{power of simple} If a free homotopy class of curves $x$ can be represented by a power of a simple curve then $\Delta(x)=0$.
\end{prop}

\section{Linked pairs}\label{review c}
In this section, we state some results from \cite{chas} which are used throughout this work.
 For convenience, we reformulate these ideas in terms of the notation of \cite{CK}.


\start{defi}{oriented}
Let $r$ be an integer greater than three and $(x_0$, $x_1$,\ldots, $x_{r-1})$ a sequence of $r$ letters in the ordered alphabet $\aq$. The sequence is \emph{ positively} (resp. \emph{negatively})
\emph{ oriented} if for some a cyclic permutation  the sequence is strictly increasing (resp. strictly decreasing).
The sign of  $(x_0,x_1,\ldots, x_{r-1})$ , denoted by $\sign(x_0,x_1,\ldots, x_{r-1})$, is 1 (resp. -1) if  $(x_0,x_1,\ldots, x_{r-1})$ is positively (resp. negatively) oriented and zero otherwise.
\end{defi}

Fix $\v=v_0\dots v_{n-1}$  a cyclically reduced word in
$\aq$ and  $\p$ and $\q$  linear subwords of $\widehat{\v}$ of equal length.


\start{defi}{linked} \cite[Definition 2.1]{chas}  The ordered pair $(\p,\q)$  is a {\em linked pair} if there exists non-negative integers $i$ and $j$ smaller that $n$ such that one of the following holds:
\begin{numlist}
\item $\p=v_{i-1}v_i$, $\q=v_{j-1}v_j$  and  $(\ov{v}_{i-1},\ov{v}_{j-1},v_i,v_j)$ is oriented. (Note that in particular $v_{i-1} \ne v_{j-1}$, $v_{i} \ne v_{j}$)
\item $\p=v_{i-r-1}\y  v_{i}$, $\q=v_{j-r-1}\y  v_{j}$ for some $r$ in  $\{1,2,\dots,n-2\}$. Moreover,
$\y=v_{i-r}v_{i-r+1}\cdots v_{i-1}=v_{j-r}v_{j-r+1}\cdots v_{j-1}$   and the sequences $(\ov{v}_{i-r-1} ,\ov{v}_{j-r-1} ,v_{i-r})$ and $(v_{i}, v_{j}
,\ov{v}_{i-1})$ have the same orientation.

\item $\p=v_{i-r-1}\y  v_{i}$, $\q=v_{j-1}\overline{\y}  v_{j+r}$ for some $r$ in  $\{1,2,\dots,n-2\}$. Moreover, $\y=v_{i-r}v_{i-r+1}\cdots v_{i-1}=\ov{v}_{j+r-1}\ov{v}_{j+r-2}\cdots \ov{v}_{j}$  and the sequences $(v_{j+r}, \ov{v}_{i-r-1},v_{i-r})$ and
$(\overline{v}_{i-1},v_{i},\ov{v}_{j-1} )$ have the same
orientation.
\end{numlist}

The \emph{ sign} of a linked pair  ${(\p,\q)}$ is given by the formula
$\sign(\p,\q)= \sign(\ov{v}_{i-1} ,v_i ,v_j) $.

The set of linked pairs of a cyclic word $\widehat{\v}$ is denoted by $\bb_1(\widehat{\v})$.

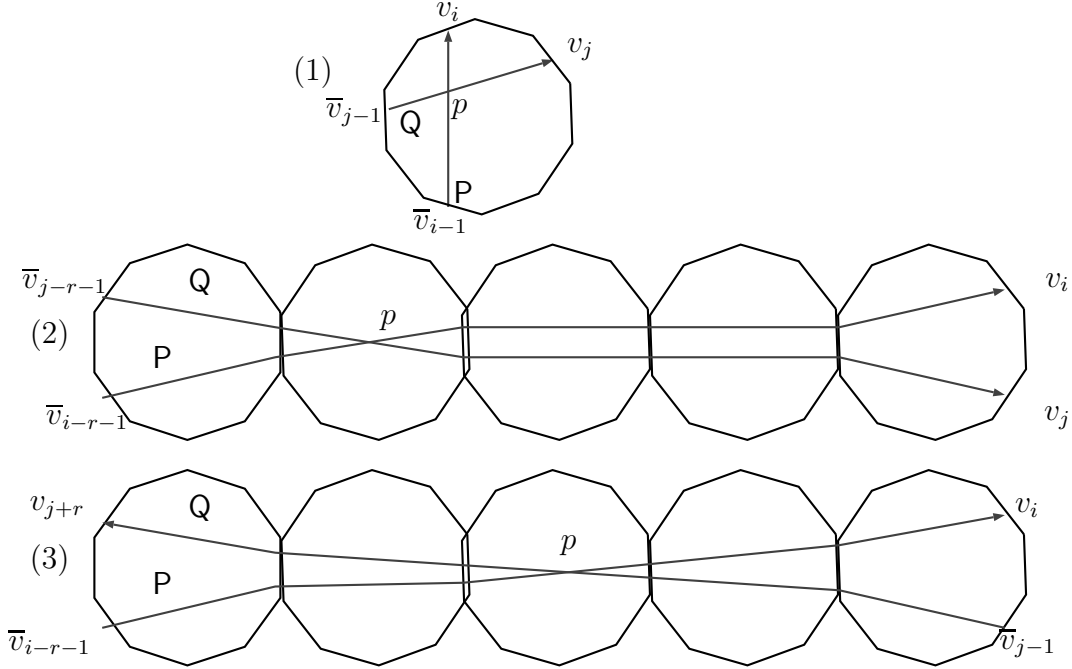
\begin{figure}[ht]
\begin{pspicture}(14,9)
\rput(6,6.4){$\p$}\rput(5.3,7.3){$\q$}
\rput(2.5,5.2){$\q$} \rput(2,4.2){$\p$}
\rput(0,-3){
\rput(2.5,5.2){$\q$} \rput(2,4.2){$\p$} }
\rput(0,3){\rput(1,0.4){$\ovv_{i-r-1}$} \rput(0.7,2.2){$\ovv_{j-r-1}$}}
\rput(0,3){\rput(13.9,0.4){$v_{j}$} \rput(13.9,2.2){$v_{i}$}}

\rput(0.5,0.4){$\ovv_{i-r-1}$} \rput(13.5,0.4){$\ovv_{j-1}$}
\rput(0.6,2.2){$v_{j+r}$} \rput(13.5,2.2){$v_{i}$}

\rput(4,8){$(1)$}\rput(0.5,4.5){$(2)$}\rput(0.5,1.5){$(3)$} \rput(-1,0){
\rput(1,6){\psline[arrowinset=0,linecolor=darkgray]{->}(5.8,0.2)(5.8,2.55)
\psline[arrowinset=0,linecolor=darkgray]{->}(5.01,1.5)(7.2,2.156)
\rput(5.7,0){$\overline{v}_{i-1}$} \rput(5.8,2.8){${v}_{i}$}
\rput(4.55,1.47){$\overline{v}_{j-1}$} \rput(7.56,2.3){$v_{j}$}

\rput{20}(6.2,1.4){\PstPolygon[unit=1.3,PolyNbSides=10]}}
\rput(2,0){ \rput(0,0){\sequence} \rput{0}(0,3){\sequence}
\psline[arrowinset=0,linecolor=darkgray]
{->}(0.2,3.66)(2.5,4.2)(5,4.6)(7.5,4.6)(10,4.6)(12.2,5.1)
\psline[arrowinset=0,linecolor=darkgray,doublesep=0.0003cm]
{->}(0.2,5)(2.5,4.6)(5,4.2)(10,4.2)(12.2,3.7)
\psline[arrowinset=0,linecolor=darkgray,doublesep=0.0003cm]
{->}(0.2,0.6)(2.5,1.15)(5,1.2)(10,1.7)(12.2,2.1)
\psline[arrowinset=0,linecolor=darkgray,doublesep=0.0003cm]
{<-}(0.2,2)(2.5,1.6)(10,1.1)(12.2,0.6)}}
\rput(5.95,7.5){$p$}\rput(5,4.7){$p$}\rput(7.4,1.7){$p$}
\end{pspicture}
\caption{Linked pairs. The label p denotes a self-intersection point.}\label{strip}
\end{figure}

\end{defi}

\start{defi}{cob12}
To each linked pair $(\p,\q)$ in $\bb_1(\widehat{\v})$ we associate two cyclic words,
$\cob_1(\p,\q)=\wh{\v}_{\hbox{\tiny $\!\!i\!,\!j$}}$ and $\cob_2(\p,\q)= \wh{\v}_{\hbox{\tiny $\!\!j\!,\!i$}}$.
\end{defi}

%
%

\start{defi}{cobr}
Let $\map{\delta}{\V}{\V\otimes \V}$  be the linear map defined by for every $\wv$ in $\V$ by the equation,
\begin{equation*}
    \cob(\wv)=\sum_{(\p,\q)\in \bb_1(\wv)}\sign(\p,\q)\,\,\cob_1(\p,\q)\otimes
    \cob_2(\p,\q).
\end{equation*}
\end{defi}

Recall that if $\Sigma$ is a surface with non-empty boundary, the set of non-trivial, free homotopy classes of closed, oriented curves on $\Sigma$ is in bijective correspondence with the set of non-empty,  cyclic reduced words on a set of free generators of the fundamental group of $\Sigma$ and their inverses.

\start{theo}{iso} (\cite[Proposition 4.1]{chas}) If $\Sigma$ be an
orientable surface with non-empty boundary then $\delta(\wv)=\Delta(\wv)$ for each $\wv$ in $\V$.
\end{theo}

The next theorem was stated as  Remark 3.10 in \cite{chas}.

\start{theo}{intersections} If $\VV$ is a
nonpower cyclic word and $\v$ a linear representative of $\VV$ then  the number of linked pairs of $\wv$ equals twice the self-intersection number of $\VV$.
\end{theo}
\begin{proof}
By \cite[Proposition 3.3]{chas}, there exists a closed curve $\alpha$ representing $\VV$ with no bigons. By \cite[Theorem 4.2]{hs} the number of self-intersection points of $\alpha$ equals the intersection number of $\VV$.
On the other hand, by \cite[Theorem 3.9]{chas} the self-intersection points of $\alpha$ are in bijective correspondence with the set of pairs of linked pairs of $\VV$ of the form $\{(\p,\q),(\q,\p)\}$. Finally, if $((\p,\q)$ is a linked pair then $\p \ne \q$. Thus the desired conclusion follows.
 \end{proof}


\section{The curve complex}\label{cc}

Let $\Sigma$ be a compact oriented surface.  By $\Sigma_{g,b}$ we denote an oriented surface with genus $g$
and $b$ boundary components. If $\Sigma$ is a surface, we denote by $\mathcal{MCG}(\Sigma)$ the \emph
{mapping class group of $\Sigma$}, that is, the set of homotopy classes of orientation preserving homeomorphisms of $
\Sigma$.  We study automorphisms of the Turaev Lie coalgebra that are related to the mapping class group.

Now we recall the curve complex,  defined by Harvey in \cite{wh}. The \emph{curve complex   $\mathrm{C}(\Sigma)$ of $\Sigma$}
 is the simplicial complex  whose vertices are  isotopy classes of unoriented  simple closed
curves on $\Sigma$ which are neither null-homotopic nor homotopic to a boundary component. If $\Sigma \ne
\Sigma_{0,4} $ and $\Sigma \ne \Sigma_{1,1}$ then a set of $k+1$ vertices of the curve complex is the $0$-
skeleton of a $k$-simplex if the corresponding  minimal intersection number of all pairs of vertices is zero, that is, if
every pair of vertices have disjoint representatives.

For $\Sigma_{0,4}$ and $\Sigma_{1,1}$ two vertices are connected by an edge when the curves they represent
have minimal intersection (2 in
the case of $\Sigma_{0,4}$ , and 1 in the case of $\Sigma_{1,1}$).  If $b \le 3$ the complex associated with $\Sigma_{0,b}$ is empty.

The following isomorphism is a theorem of Ivanov
\cite{iv} for the case of genus at least two.  Korkmaz \cite{ko} proved the result for genus at most one and Luo \cite{fl} gave
another proof  that covers all possible genera. The mapping class group of a surface $\Sigma$ is denoted by $\mathcal{MCG}(\Sigma)$. Our statement below is based on the formulation of Minsky \cite{minsky}.

\start{theo}{ikl}(Ivanov-Korkmaz-Luo) The natural map $\map{h}{\mathrm{\mathcal{MCG}}(\Sigma)}{\mathrm{Aut C}(\Sigma)}$ is an isomorphism in all cases except for $\Sigma_{1,2}$ where it is injective with index $2$ image.
 \end{theo}

\end{document}